\documentclass[10pt]{amsart}
\usepackage{amsthm,amsmath,amssymb,amsfonts}
\usepackage[american]{babel}
\usepackage{cite}
\usepackage{url}

\newtheorem{theorem}{Theorem}[section] 
\newtheorem{lemma}[theorem]{Lemma}
\newtheorem{proposition}[theorem]{Proposition}

\newtheorem{claim}[theorem]{Claim}

\theoremstyle{definition}

\newtheorem{example}[theorem]{Example}

\theoremstyle{remark}
\newtheorem{remark}[theorem]{Remark}

\numberwithin{equation}{section}

\newcommand{\abs}[1]{\lvert#1\rvert}

\DeclareMathOperator{\supp}{supp}

\begin{document}


\newcommand{\nt}{\trianglelefteq}
\newcommand{\aut}[1]{\text{\rm Aut}(#1)}
\newcommand{\out}[1]{\text{\rm Out}(#1)}
\newcommand{\autc}[1]{\mbox{\rm Aut}_{\mbox{\rm\scriptsize c}}(#1)}
\newcommand{\outc}[1]{\mbox{\rm Out}_{\mbox{\rm\scriptsize c}}(#1)}
\newcommand{\inn}[1]{\mbox{\rm Inn}(#1)}
\newcommand{\cen}[2]{\mbox{\rm C}_{#1}(#2)}
\newcommand{\nor}[2]{\mbox{\rm N}_{#1}(#2)}
\newcommand{\zen}[1]{\mbox{\rm Z}(#1)}
\newcommand{\cyk}[1]{\text{\rm C}_{#1}}
\newcommand{\lara}[1]{\langle{#1}\rangle}
\newcommand{\ZZ}{\mbox{$\mathbb{Z}$}}
\newcommand{\NN}{\mbox{$\mathbb{N}$}}
\newcommand{\Zn}[2]{\mbox{\rm Z}_{#1}(#2)}
\newcommand{\R}[1]{\mbox{\rm R}(#1)}
\newcommand{\U}{\mathcal{U}}


\title[The normalizer property for Blackburn groups]
{Class-preserving automorphisms and the normalizer
property for Blackburn groups}

\author{Martin Hertweck}
\address{Universit\"at Stuttgart, Fachbereich Mathematik,
Institut f\"ur Geometrie und Topo\-logie,
Pfaffenwald\-ring 57, 70550 Stuttgart, Germany}
\email{hertweck@mathematik.uni-stuttgart.de}

\author{Eric Jespers }
\address{Department of Mathematics,
Vrije Universiteit Brussel,  Pleinlaan 2, 1050
Brussel, Belgium} \email{efjesper@vub.ac.be}
\thanks{The second author is partially supported by
the Onderzoeksraad of Vrije Universiteit Brussel,
Fonds voor Wetenschappelijk Onderzoek (Flanders) and
the Flemish-Polish bilateral agreement
BIL2005/VUB/06.}


\date{\today}


\begin{abstract}
For a group $G$, let $\U$ be the group of units of
the integral group ring $\ZZ G$. The group $G$
is said to have the normalizer property if
$\nor{\U}{G}=\zen{\U}G$. It is shown that Blackburn
groups have the normalizer property. These are the
groups which have non-normal finite subgroups, with
the intersection of all of them being nontrivial.
Groups $G$ for which class-preserving automorphisms
are inner automorphisms, $\outc{G}=1$, have the
normalizer property. Recently, Herman and Li have
shown that $\outc{G}=1$ for a finite Blackburn group
$G$. We show that $\outc{G}=1$ for the members $G$
of a few classes of metabelian groups, from which
the Herman--Li result follows.

Together with recent work of Hertweck, Iwaki,
Jespers and Juriaans, our main result implies that,
for an arbitrary group $G$, the group
$\Zn{\infty}{\U}$ of hypercentral units of $\U$
is contained in $\zen{\U}G$.
\end{abstract}

\maketitle

\section{Introduction}

A group $G$ is said to have the {\em normalizer
property} if $\nor{\U}{G}=\zen{\U}G$, where $\U$
denotes the group of units of the integral
group ring $\ZZ G$, $\zen{\U}$ denotes the center of
$\U$ and $\nor{\U}{G}$ is the normalizer of $G$ in
$\U$. Informally speaking, this means that $G$ is
normalized only by units of $\U$ which obviously do
so. Sometimes this definition is slightly tightened
to incorporate $G$-adapted coefficient rings: these
are integral domains of characteristic zero in which
a rational prime $p$ is not invertible whenever $G$
has an element of order $p$. Throughout, we let $R$
denote a $G$-adapted ring, write $\U$ for the group
of units of $RG$, and if $\nor{\U}{G}=\zen{\U}G$, we
say that $G$ has the normalizer property. For
classes of groups which are known to have the
normalizer property, we refer the reader to
\cite{HIJJ:06} (and references therein).

The normalizer $\nor{\U}{G}$ appears, among other
things, in the study of the group $\Zn{\infty}{\U}$
of hypercentral units of $\U$, defined by
$\Zn{\infty}{\U}=\bigcup_{n=1}^{\infty}\Zn{n}{\U}$,
where $\Zn{n}{\U}$ denotes the $n$-th term of the
upper central series of $\U$. The question whether
in general $\Zn{\infty}{\U}\leq\nor{\U}{G}$ was
finally answered in the affirmative \cite{HIJJ:06}.
As a consequence of the present work we obtain that
$\Zn{\infty}{\U}\leq\zen{\U}G$. We remark that this
is known to hold for $G$ torsion-free
\cite[Corollary~4.3]{HIJJ:06}, and for $G$ whose
finite subgroups are all normal
\cite[Corollary~4.12]{HIJJ:06}. If $G$ has finite
non-normal subgroups, then, following Blackburn
\cite{Bl:66}, we denote by $\R{G}$ the intersection
of all of them and say that $\R{G}$ is defined. Now,
if $\R{G}$ is defined and $\R{G}=1$, then
$\Zn{\infty}{\U}=\zen{\U}$ by \cite[Propositions~4.1
and 4.5]{HIJJ:06}. It remains to consider the {\em
Blackburn groups}, i.e., the groups $G$ for which
$\R{G}$ is defined and nontrivial. Here, it is shown
that Blackburn groups have the normalizer property.

We write $\autc{G}$ for the group of
class-preserving automorphisms of $G$, and set
$\outc{G}=\autc{G}/\inn{G}$. Any $u$ in
$\nor{\U}{G}$ gives rise---via conjugation---to a
class-preserving automorphism of $G$ (see
\cite[Proposition~2.6]{HIJJ:06}). Thus $\outc{G}=1$
is a sufficient condition for $G$ to have the
normalizer property. In praxis, this has proven to
be a valuable criterion. We would have to point out,
though, that we do not know whether class-preserving
automorphisms of infinite Blackburn groups are inner
automorphisms.

The finite Blackburn groups were classified in
\cite{Bl:66}. Li, Sehgal and Parmenter
\cite[Theorem~1]{LiSePa:99} have shown by inspection
that the finite Blackburn groups have the normalizer
property. In fact, $\outc{G}=1$ for a finite
Blackburn group $G$. Herman and Li noted that the
proof given in \cite[22.4~Proposition]{He:04a}
contains a gap, and completed it in \cite{HeLi:06}
(Lemma~\ref{again2} below also fixes the gap).

In Section~\ref{Sec2:CPAut}, we record that
class-preserving automorphisms of finite
abelian-by-cyclic groups are inner automorphisms
(see Proposition~\ref{PrAC}). Together with an ad
hoc observation (Lemma~\ref{again2}), the
last-mentioned result again follows.

In Section~\ref{Sec3:BGroupsNP}, the results from
\cite[\S3]{HIJJ:06} are applied to show that
infinite Blackburn groups have the normalizer
property (Theorem~\ref{peve}). In doing so, an
essential feature will be that the structure of the
finite normal subgroups of Blackburn groups is known
since these are either Dedekind groups or again
Blackburn groups.

Finally, let us glance over Blackburn's
classification \cite{Bl:66}. By
\cite[Theorem~1]{Bl:66}, a finite $p$-group $G$
which is a Blackburn group is a $2$-group, and
either a Q-group or of the form
$\text{Q}_{8}\times\text{C}_{4}\times\text{E}_{2}$
or
$\text{Q}_{8}\times\text{Q}_{8}\times\text{E}_{2}$,
where $\text{Q}_{8}$ is the quaternion group,
$\text{E}_{2}$ is an elementary abelian $2$-group
and $\text{C}_{n}$ denotes the cyclic group of order
$n$. (A group $G$ is called a {\em Q-group} if $G$
has an abelian subgroup $A$ of index $2$ which is
not elementary abelian, and $G=\lara{A,b}$ for some
$b\in G$ of order $4$ with $x^{b}=x^{-1}$ for all
$x\in A$.) A finite Blackburn group which is not a
$p$-group belongs to one of five classes which are
listed in \cite[Theorem~2]{Bl:66}. In particular,
such a group is the direct product of groups which
are abelian or possess an abelian normal subgroup
with factor group cyclic or isomorphic to Klein's
four group.

\section{On class-preserving automorphisms}\label{Sec2:CPAut}

In this section, only finite groups are considered.
We present classes of metabelian groups for
which class-preserving automorphisms are inner
automorphisms, and conclude by showing
that the finite Blackburn groups belong to them.

We first prove that class-preserving automorphisms
of finite abelian-by-cyclic groups are inner
automorphisms. For the class of metacyclic groups,
this is an elementary exercise. Also, the following
example, used in the proof of Lemma~\ref{again2}, is
easily verified (see \cite[Example~3.1]{HIJJ:06}).
\begin{example}\label{again1}
A group $G$ with an abelian subgroup of index $2$
satisfies $\outc{G}=1$.
\end{example}

Our proof will be based on the following lemma,
which is of interest in its own right.

\begin{lemma}\label{ispower}
Let $G$ be a finite abelian $p$-group, and let
$\alpha$ and $\beta$ be automorphisms of $G$ of
$p$-power order. Assume that
$\alpha\beta=\beta\alpha$ and that for each $g\in
G$, there is $n\in\mathbb{N}$ such that
$g\beta=g\alpha^{n}$. Then $\beta$ is a power of
$\alpha$.
\end{lemma}
\begin{proof}
Assume that $G$ is a counterexample, with the order
of the semidirect product $G\langle \alpha\rangle$
being minimal. Let $Z$ be a central subgroup of
order $p$ in $G\langle \alpha\rangle$ which is
contained in $G$. Then $\beta$ centralizes $Z$. Thus
$\alpha$ and $\beta$ induce automorphisms
$\bar{\alpha}$ and $\bar{\beta}$ of $\bar{G}=G/Z$,
and $\bar{\beta}$ is a power of $\bar{\alpha}$ by
minimality of $G\langle \alpha\rangle$. Hence we can
assume that $\bar{\beta}$ is the identity. Then the
map $G\rightarrow Z$, defined by $g\mapsto
g^{-1}(g\beta)$, is a surjective homomorphism, with
kernel $K$ of index $p$. Clearly $K=\cen{G}{\beta}$
and $\alpha$ fixes $K$ since $\alpha$ and $\beta$
commute. For all $g\in G$, there is
$n(g)\in\mathbb{N}$ such that
$g\beta=g\alpha^{n(g)}$. Choose $h\in G\setminus K$
such that $\alpha^{n(h)}$ is of maximal order among
the $\alpha^{n(g)}$, $g\in G\setminus K$. By
minimality of $G\langle \alpha\rangle$, we have
$\langle\alpha\rangle= \langle\alpha^{n(h)}\rangle$.
Thus we can assume that $G=\langle K,h\rangle$ and
$h\beta=h\alpha$. Note that $h\beta=zh$ for a
central element $z$ of order $p$ in $Z$.

For all $k\in K$, we have
$kzh=(k\beta)(h\beta)=(kh)\beta=(kh)\alpha^{n(kh)}=
k\alpha^{n(kh)}z^{n(kh)}h$. So $\alpha^{n(kh)}$ is
not the identity and $k\alpha^{n(kh)}\in kZ$. Thus
if $\alpha^{q}$ is a power of $\alpha$ having order
$p$, then $\alpha^{q}$ induces the identity on
$G/Z$. As above, it follows that $H:=\text{\rm
C}_{G}(\alpha^{q})$ is an $\alpha$-invariant
subgroup of $G$ of index $p$. Also, $H$ is fixed by
$\beta$ since $\beta$ commutes with $\alpha$. Thus
by minimality of $G\langle \alpha\rangle$, the
automorphism $\beta$ agrees on $H$ with some power
$\alpha^{l}$ of $\alpha$. Suppose that $h\not\in H$.
Then $h\neq h\alpha^{q}=z^{q}h$, so $p\nmid q$ and
$\alpha$ has order $p$. Consequently $H=\text{\rm
C}_{G}(\alpha)\leq\text{\rm C}_{G}(\beta)=K$, so
$G=\lara{h,H}$ and $\beta=\alpha$, contradicting our
assumption that $G$ is a counterexample. Thus $h\in
H$. It follows that $zh=h\beta=h\alpha^{l}$, and 
$\langle\alpha^{l}\rangle=\langle\alpha\rangle$ 
since $h$ is a fixed point under $\alpha^{p}$.

Since $\beta$ induces the identity on $K$, it
follows that $\alpha^{l}$, and hence also $\alpha$,
induces the identity on $H\cap K$. 
Take any $k\in K\setminus H$. Since $H\cap K$ is
of index $p^{2}$ in $G$, it follows that $G=\langle
h,k,H\cap K\rangle$. Since $\alpha$ induces the
identity on $G/H$ (a quotient of order $p$) and
$K$ is $\alpha$-invariant, 
we have $k\alpha=xk$ for some $x\in H\cap K$, 
and $x\neq 1$ since $\alpha\neq\beta$. Let
$x^{m}$ be a power of $x$ having order $p$. Suppose
that $\langle x^{m}\rangle=\langle z\rangle$, i.e.,
$x^{mt}=z^{-1}$ for some $t\in\NN$. Then
$(hk^{mt})\alpha=hk^{mt}$ 
(this is the crucial place where we need that
$G$ is abelian) but $(hk^{mt})\beta\neq hk^{mt}$, a 
contradiction. Thus 
$\langle x\rangle\cap \langle z\rangle=1$.

Since $x\alpha=x=x\beta$, the automorphisms $\alpha$
and $\beta$ induce automorphisms of $G/\langle
x^{m}\rangle$, and by minimality of $G\langle
\alpha\rangle$, the automorphism induced by $\beta$
equals the automorphism induced by some power
$\alpha^{j}$ of $\alpha$. Then
$z^{j-1}=h\alpha^{j}(h\beta)^{-1}\in\langle
x^{m}\rangle$ and thus $j=1+ps$ for some
$s\in\mathbb{Z}$. Furthermore,
$x^{j}=x^{j}kk^{-1}=k\alpha^{j}(k\beta)^{-1}\in\langle
x^{m}\rangle$. So $\langle x\rangle=\langle
x^{j}\rangle\leq\langle x^{m}\rangle$ and $x$ has
order $p$. Because $\langle x \rangle \cap \langle z
\rangle =1$, it follows that
$(hk)\alpha^{i}=z^{i}hx^{i}k\neq zhk=(hk)\beta$ for 
all $i\in\mathbb{N}$, and we
have reached a final contradiction.
\end{proof}

In the above lemma, the assumption that $G$ is 
abelian is necessary. In fact, the proof lends 
itself to the construction of adequate examples.
Such examples yield $p$-groups with non-inner,
class-preserving automorphisms.

\begin{example}
For an odd prime $p$, we construct a group $G$ of order
$p^{p+2}$ which has automorphisms $\alpha$ and $\beta$ 
satisfying the assumptions, but not the conclusion
of the above lemma
(for $p=2$, see \cite[14.3 Remark]{He:04a}).
First, we have to exhibit a certain kind of group $K$
which could not appear in the above (abelian group) setting. Let
\renewcommand{\arraystretch}{-2}
\[ \kappa = \left( \!
\begin{array}{c@{\hspace*{2.5pt}}c@{\hspace*{2.5pt}}c
@{\hspace*{2.5pt}}c@{\hspace*{2.5pt}}c@{\hspace*{2.5pt}}c}
1 & p & 0 & \cdots & \cdots & 0 \vspace*{-2pt} \\
\raisebox{2pt}{0} & \raisebox{2pt}{1} & \raisebox{2pt}{1} 
& \ddots & & \vdots \vspace*{-3.5pt} \\
\vdots & \ddots & \ddots & \ddots & \ddots & 
\vdots \vspace*{-3.5pt} \\
\vdots & & \ddots & \ddots & 1 & 0 \vspace*{-3.5pt} \\
0 & & & \ddots & 1 & 1 \vspace*{5pt} \\
-1 & 0 & \cdots & \cdots & 0 & 1
\end{array} 
\! \right) \in \text{\rm GL}_{p-1}(\ZZ/p^{2}\ZZ), \]
\renewcommand{\arraystretch}{1}
an invertible $(p-1)\times(p-1)$ matrix (of order $p^{2}$)
over $\ZZ/p^{2}\ZZ$, which acts on
\[ A=\ZZ/p^{2}\ZZ\oplus\underset{p-2\text{ times}}
{\underbrace{p\ZZ/p^{2}\ZZ\oplus\ldots\oplus p\ZZ/p^{2}\ZZ}} \cong 
\cyk{p^{2}}\times\cyk{p}\times\ldots\times\cyk{p} \]
as an automorphism of order $p$, and has the property that 
$a(1+\kappa+\ldots+\kappa^{p-1})=0$ for all $a\in A$
(cf.\ \cite[III.10.15]{Hupp:67}).
For completeness, we give the argument: The integer 
polynomials $f(t)=t^{p-2}+2t^{p-3}+\ldots+(p-2)t+(p-1)$ and
$c(t)=t^{p-1}+t^{p-2}+\ldots+1$ satisfy $(t-1)f(t)=c(t)-p$.
So, modulo $p$, we have
$(t-1)^{2}f(t)\equiv t^{p}-1\equiv(t-1)^{p}$ and
$f(t)\equiv (t-1)^{p-2}$. Let $a\in A$.
Then $pa(\kappa-1)=0$, and so
$a(\kappa-1)f(\kappa)=a(\kappa-1)^{p-1}$ by the 
established congruence. Note that $(\kappa-1)^{p-1}$ is $-p$
times the identity matrix. Altogether, we obtain 
$a(1+\kappa+\ldots+\kappa^{p-1})=ac(\kappa)=
a(\kappa-1)f(\kappa)+pa=0$. From that, it is obvious that 
$\kappa$ acts as an automorphism of order $p$ on $A$.
We are ready to present the group construction.
Let $k$ be an (abstract) generator of a cyclic group of order
$p$, and form the semidirect product $K=A\langle k\rangle$,
with $k$ acting on $A$ via $\kappa$. Let $x$ be a generator 
of the first factor of $A$, and set $z=x^{p}$, a central
element of order $p$ in $K$. By the mentioned property of 
$\kappa$ (with $a=x$), the element $xk$ has order $p$, so that 
an automorphism $\alpha$ of $K$ is defined by prescribing that
$\alpha$ fixes $A$ point-wise, and $k\alpha=xk$.

Now let $h$ be another (abstract) generator of a cyclic group
of order $p$, and set $G=K\times\langle h\rangle$.
Extend $\alpha$ to an automorphism of $G$ by setting
$h\alpha=zh$. Let $\beta$ be the automorphism of $G$ 
which fixes $K$ point-wise, and sends $h$ to $zh$.
Certainly, $\alpha$ has order $p^{2}$, $\beta$ has order $p$,
and both automorphisms commute.
Note that $H:=\cen{G}{\alpha^{p}}=A\times\langle h\rangle$.
Let $g\in G\setminus K$. If $g\in H$, then $g\beta=g\alpha$.
Otherwise $g\in Ak^{i}h^{j}$ with $1\leq i,j\leq p-1$, and
$g\beta=g\alpha^{ps}$ for $s\in\NN$ with $is\equiv j\mod{p}$.
However, $\beta$ is not a power of $\alpha$.

It readily follows that the automorphism of
the semidirect product $G\langle\alpha\rangle$
which fixes $\alpha$ and agrees on $G$ with $\beta$
is a non-inner, class-preserving automorphism.
\end{example}

We shall need the following well known lemma for
which we did not found a suitable reference in
the literature.
\begin{lemma}\label{HN}
If $N$ is a $p$-group on which an abelian
$p^{\prime}$-group $H$ acts then
$\cen{H}{n_{0}}=\cen{H}{N}$ for some $n_{0}\in N$.
\end{lemma}
\begin{proof}
Let $\bar{N}$ be the quotient of $N$ by its Frattini
subgroup. Then $\cen{H}{\bar{N}}=\cen{H}{N}$ (see
\cite[(24.1)]{Asch:86}). So  it suffices to find
$n_{0}\in N$ with
$\cen{H}{\bar{n}_{0}}=\cen{H}{\bar{N}}$. Thus we can
assume that $N$ is elementary abelian. Form the
semidirect product $X=HN$. By Maschke's theorem,
$N=M_{1}\times\ldots\times M_{r}$ with minimal
normal subgroups $M_{i}$ of $X$. Choose $m_{i}\in
M_{i}\setminus\{1\}$ for each index $i$. Note that
for $h\in H$, each subgroup $\cen{M_{i}}{h}$ is
normal in $X$ since $H$ is abelian. So if $h\in H$
with $m_{i}^{h}=m_{i}$ for some $i$ then
$\cen{M_{i}}{h}\neq 1$ and $\cen{M_{i}}{h}=M_{i}$ as
$M_{i}$ is a minimal normal subgroup. In other
words, $\cen{H}{m_{i}}=\cen{H}{M_{i}}$ for all $i$.
Hence for $h\in\cen{H}{m_{1}\cdots m_{r}}$ we have
$h\in\bigcap_{i=1}^{r}\cen{H}{m_{i}}=
\bigcap_{i=1}^{r}\cen{H}{M_{i}}=\cen{H}{N}$, and the
lemma holds with $n_{0}=m_{1}\cdots m_{r}$.
\end{proof}

For the reader's convenience, we recall the
following easy but useful observation (see
\cite[Remark~1]{He:01b}).
\begin{remark}\label{rem1}
Let $\sigma$ be a class-preserving automorphism of
$G$ of order a power of $p$, and let
$\gamma\in\inn{G}$. Then $(\gamma\sigma)^{r}$, where
$r$ denotes the $p^{\prime}$-part of the order of
$\gamma\sigma$, is again a class-preserving
automorphism of $G$ of order a power of $p$ which is
a non-inner automorphism if and only if $\sigma$ is
a non-inner automorphism (since
$\lara{\sigma^{r}}=\lara{\sigma}$ and
$(\gamma\sigma)^{r}$ and $\sigma^{r}$ differ only by
an inner automorphism).

Thus if one wishes to show that such a $\sigma$ is
an inner automorphism, and $U$ is a subset of $G$
which is conjugate to $U\sigma$ in $G$, one can
assume in addition that $U=U\sigma$.
\end{remark}

Now we are in a position to prove:
\begin{proposition}\label{cyclicp}
Let $G$ be a finite group having an abelian normal
subgroup $A$ such that the quotient $G/A$ has a
normal cyclic Sylow $p$-subgroup, for some prime
$p$. Then each class-preserving automorphism of $G$
of $p$-power order is an inner automorphism.
\end{proposition}
\begin{proof}
Let $\sigma$ be a class-preserving automorphism of
$G$ of $p$-power order; we have to show that
$\sigma$ is an inner automorphism. Let $P$ be a
Sylow $p$-subgroup of $G$. By
\cite[Corollary~5]{He:01b} (applied with $N=PA$) we
can assume that $G/A$ is a $p$-group. Then we can
choose $x\in P$ such that $G=\langle x,A\rangle$. By
Sylow's theorem, and Remark~\ref{rem1}, we can
assume that $P\sigma=P$. Set $S=P\cap A=\text{\rm
O}_{p}(A)$ and $T=\text{\rm O}_{p^{\prime}}(A)$, so
that $P=\langle x,S\rangle$ and $A=S\times T$.

Let $\gamma\in\aut{G}$ be the inner automorphism
given by conjugation with $x$ and set $H=\langle
\sigma|_{T},\gamma|_{T}\rangle\leq\text{\rm
Aut}(T)$. Note that for each $t\in T$ we have
$t\sigma|_{T}=t(\gamma|_{T})^{n}$ for some
$n\in\NN$, and that $x\sigma\in xA\leq x\cen{G}{T}$.
Thus $H$ is an abelian $p$-group, $\text{\rm
C}_{H}(t_{0})=\text{\rm C}_{H}(T)=1$ for some
$t_{0}\in T$ by Lemma~\ref{HN}, and $\sigma|_{T}$ is
a power of $\gamma|_{T}$. So Remark~\ref{rem1}
allows us to assume that $t\sigma=t$ for all $t\in
T$.

Let $y$ be a generator of $\text{\rm C}_{\langle
x\rangle}(t_{0})$, and let $\delta$ be the inner
automorphism given by conjugation with $y$. For each
$s\in S$ there is $n(s)\in\mathbb{N}$ such that
$(s\sigma)
t_{0}=(st_{0})\sigma=(st_{0})\gamma^{n(s)}=
(s\gamma^{n(s)})(t_{0}\gamma^{n(s)})$, meaning that
$\gamma^{n(s)}\in\lara{\delta}$ and
$s\sigma=s\delta^{m(s)}$ for some
$m(s)\in\mathbb{N}$. As before one sees that
$\sigma|_{S}$ commutes with $\delta$. Thus
$\sigma|_{S}$ is a power of $\delta$ by
Lemma~\ref{ispower}, and we can modify $\sigma$
according to Lemma~\ref{HN} such that the new
$\sigma$ fixes $A$ element-wise. Clearly
$x\sigma=x^{s}$ for some $s\in S$ since $x\sigma\in
P$, and then $\sigma$ is the inner automorphism
given by conjugation with $s$.
\end{proof}

From that we obtain at once:
\begin{proposition}\label{PrAC}
Let $G$ be a finite group having an abelian normal
subgroup $A$ with cyclic quotient $G/A$. Then
class-preserving automorphisms of $G$ are inner
automorphisms.
\end{proposition}
\begin{proof}
Let $\sigma$ be a class-preserving automorphism of
$G$ of $p$-power order, for some prime $p$; we have
to show that $\sigma$ is an inner automorphism. By
\cite[Corollary~5]{He:01b} we can assume that $p$
divides the order of $G/A$, and then
Proposition~\ref{cyclicp} applies.
\end{proof}

In preparation for the final result we record:
\begin{lemma}\label{again2}
Suppose that a finite group $G$ is a semidirect
product of an abelian group $A$ and a generalized
quaternion group $Q$, and that each element of $Q$
acts on $A$ by raising each element to some fixed
power. Suppose further that a Sylow $2$-subgroup of
$G$ has an abelian subgroup of index $2$. Then
$\outc{G}=1$.
\end{lemma}
\begin{proof}
Let $\sigma\in\autc{G}$. We shall show that
$\sigma\in\inn{G}$ by induction on the order of $G$.
By Sylow's theorem, we can assume that $S\sigma=S$
for a Sylow $2$-subgroup $S$ of $G$ containing $Q$.
Then $\sigma$ induces a class-preserving
automorphism of $S$ since $S$ can be viewed as a
homomorphic image of $G$. Thus, by
Example~\ref{again1}, we can assume that $\sigma$
fixes $S$ element-wise. Since the conjugation action
of $Q$ on $A$ leaves each cyclic subgroup invariant,
and the automorphism group of a cyclic group is
abelian, the commutator subgroup $Q^{\prime}$
centralizes $A$.

Suppose that $A=M\times N$ with nontrivial normal
subgroups $M$ and $N$ of $G$. Then we can assume
inductively that $\sigma$ induces an inner
automorphism of $G/N$, say conjugation with $Ng$
(some $g\in G$). Write $g=xa$ with $x\in Q$ and
$a\in A$. Then for any $y\in Q$, since $Q\leq S$ and
$\sigma$ fixes $S$ element-wise, we have
$Ay=A(y\sigma)=Ay^{xa}=Ay^{x}$, i.e.,
$y^{x}y^{-1}\in Q\cap A=1$. So $x\in\zen{Q}\leq
Q^{\prime}\leq\cen{G}{A}$. It follows that for $m\in
M$, we have $N(m\sigma)=Nm^{xa}=Nm^{a}=Nm$. Hence
$\sigma$ induces the identity on $G/N$. Likewise, we
obtain that $\sigma$ induces the identity on $G/M$,
so $\sigma$ is the identity.

Thus we can assume that $A$ is cyclic $p$-group, and
by Example~\ref{again1}, we can assume that $p$ is
an odd prime. Let $a$ be a generator of $A$. Note
that $\aut{A}$ is cyclic, so we can choose, since
$Q/Q^{\prime}$ is a Klein's four-group, an element
$y\in Q\setminus Q^{\prime}$ which commutes with
$A$. Then $(ay)\sigma=(ay)^{x}$ for some $x\in Q$,
that is, $a\sigma=a^{x}$ and $y=y\sigma=y^{x}$. The
latter implies that $\lara{x,y}$ is a proper
subgroup of $Q$, so $x\in Q^{\prime}\cup
Q^{\prime}y$ and $a\sigma=a^{x}=a$. Again, we have
shown that $\sigma$ is the identity, and the proof
is complete.
\end{proof}

Now we can go through Blackburn's list to obtain:
\begin{proposition}\label{hypex}
A class-preserving automorphism of a finite
Blackburn group is an inner automorphism.
\end{proposition}
\begin{proof}
Note that a class-preserving automorphism stabilizes
every normal subgroup, and that $\autc{-}$ commutes
with taking direct products.

Let $G$ be a finite Blackburn group. If $G$ is a
$p$-group, then $G$ contains an abelian subgroup of
index $2$. Hence, the statement follows from
Example~\ref{again1} (it is, of course, possible to
give a more direct proof).

Suppose that $G$ is not of prime-power order. Then
$G$ belongs to one of the classes (a)--(e) described
in \cite[Theorem~2]{Bl:66}. If $G$ is of type (a) or
(d), then $G$ is abelian-by-cyclic and
Proposition~\ref{PrAC} applies (of course, one can
argue more directly). The groups of type (b) or (e)
are now recognized as direct products of groups
which have only inner class-preserving
automorphisms. Thus it remains to deal with $G$ of
type (c), but such a group $G$ is of the form
described in Lemma~\ref{again2} (note that the
groups of type (a) have an abelian Sylow
$p$-subgroup), so we are done.
\end{proof}

\section{Blackburn groups have the normalizer property}\label{Sec3:BGroupsNP}

From now on, $G$ will denote a Blackburn group, so
$\R{G}$ is defined and ${\text{\rm R}}(G)\neq 1$. It
is easily seen that ${\text{\rm R}}(G)$ is a cyclic
$p$-group for some prime $p$ (see
\cite[22.1~Lemma]{He:04a} for a proof). Note that if
a finite subgroup $S$ of $G$ is not a Dedekind
group, then ${\text{\rm R}}(G)\leq{\text{\rm
R}}(S)$; in particular, ${\text{\rm R}}(S)\neq 1$.
The following observation was also stated in
\cite[Lemma~4.10]{HIJJ:06}. For completeness sake we
include the proof.

\begin{lemma}\label{qdifp}
Suppose that $G$ has a finite non-normal subgroup,
and that ${\text{\rm R}}(G)$ is a nontrivial cyclic
$p$-group. Let $q$ be a prime different from $p$.
Then the set $T_{q}$ of $q$-elements of $G$ forms a
normal $q$-subgroup of $G$, all of whose subgroups
are normal in $G$. (In particular, $T_{q}$ for odd
$q$ is an abelian $q$-group.)
\end{lemma}
\begin{proof}
Let $x$ and $y$ be $q$-elements in $G$. Then
$\lara{x}$ is normal in $G$ since otherwise we would
have ${\text{\rm R}}(G)=1$. Thus $xy$ is also a
$q$-element, and it follows that $T_{q}$ is a group.
If a subgroup of $T_{q}$ is non-normal in $G$, then
$T_{q}$ also contains a finite cyclic subgroup which
is non-normal in $G$, which again contradicts
${\text{\rm R}}(G)\neq 1$. That $T_{q}$ for odd $q$
is abelian holds by the classification of Dedekind
groups.
\end{proof}

We shall use the following information on finite
normal subgroups of $G$. Rather than deducing it
from Blackburn's list \cite{Bl:66}, we extract the
proof directly from \mbox{Blackburn}'s work.
\begin{claim}\label{fnsgp}
Suppose that $G$ has a finite non-normal subgroup,
and that ${\text{\rm R}}(G)$ is a nontrivial cyclic
$p$-group. Let $N$ be a finite normal subgroup of
$G$. Then $N$ has a normal $p$-complement $H$ which
is a Dedekind group. Let $S$ be a Sylow $p$-subgroup
of $N$. Then one of the following holds:
\begin{enumerate}
\item[\mbox{(a)}]
$N$ is nilpotent;
\item[\mbox{(b)}]
$\zen{S}\leq\text{\rm O}_{p}(N)$, i.e.,
$[\zen{S},H]=1$;
\item[\mbox{(c)}]
$S$ is abelian, and $S=\lara{x}\times T$, where
$x\in S$ is such that $[x,H]\neq 1$, and
$T\leq\text{\rm O}_{p}(N)$. The order of
$\cen{\lara{x}}{H}$ is the exponent of $\text{\rm
O}_{p}(N)$. If some $g$ in $G$ centralizes
$\cen{\lara{x}}{H}$, then $g$ centralizes $T$.
\end{enumerate}
\end{claim}
\begin{proof}
For a prime $q$ different from $p$, a Sylow
$q$-subgroup of $N$ is a normal subgroup of $N$ by
Lemma~\ref{qdifp}. So $N$ has a normal
$p$-complement $H$. Furthermore, all subgroups of
$H$ are normal in $G$  and thus  $H$ is a Dedekind
group.

If $N$ is a Dedekind group, then (a) holds. Hence we
can assume that ${\text{\rm R}}(N)\neq 1$, i.e., $N$
is a finite Blackburn group with $R(G)\leq R(N)$ and
$R(N)$ is a cyclic $p$-group. Let $S$ be a Sylow
$p$-subgroup of $N$.

Assume $S$ is not a Dedekind group. Then, $R(G)\leq
R(S)$. Thus, $S$ is a  Blackburn group and a
$p$-group. So $p=2$. If $S$ is of the form
$\text{Q}_{8}\times\text{C}_{4}\times\text{E}_{2}$
or
$\text{Q}_{8}\times\text{Q}_{8}\times\text{E}_{2}$
(where $\text{E}_{2}$ is an elementary abelian
$2$-group), then $S$ is generated by certain
subgroups, none of which contains ${\text{\rm
R}}(S)$. Since $|R(S)|=2$, we get that $R(S)=R(G)$.
Hence, the mentioned subgroups generating $S$ are
all normal in $G$. So $S$ is normal in $G$, and (a)
holds. If $S$ is a Q-group (possibly a non-abelian
Dedekind group), then $\zen{S}$ is elementary
abelian, and (b) holds since elements of $S$ of
order $2$ are contained in the center of $N$ as
${\text{\rm R}}(N)\neq 1$. If $S$ is a non-abelian
Dedekind group then the last argument also yields
that (b) holds. Hence, we are left to deal with $S$
an abelian group.

We can assume that $[S,H]\neq 1$ since otherwise (a)
holds. Set $C=\text{\rm O}_{p}(N)=\cen{S}{H}$. We
will follow \cite[p.~35]{Bl:66}. Note that if $s$ is
any element of $S$ not lying in $C$, then $\lara{s}$
is not normal in $N$ and ${\text{\rm
R}}(N)\leq\lara{s}$. Hence, ${\text{\rm
R}}(N)\leq\lara{s^{p}}$ since ${\text{\rm R}}(N)\nt
N$.

Let $x,y\in S$ with $x\not\in C$ and the order of
$x$ larger than or equal to the order of $y$. By the
basis theorem for abelian groups,
$\lara{x,y}=\lara{x}\times\lara{z}$ for some $z\in
S$, and since ${\text{\rm R}}(N)\leq\lara{x}$, we
have $z\in C$. Thus $S/C$ is a cyclic group, say of
order $p^{r}$. Suppose $S=\lara{C,x}$ and that $x$
is of order $p^{m}$. Then $m>r$ since ${\text{\rm
R}}(N)\leq\lara{x}$, and $x^{p^{m-1}}\in{\text{\rm
R}}(N)$.

We prove next that the exponent of $C$ is $p^{m-r}$.
If this is not so, then $C$ contains an element $y$
of order $p^{m-r+1}$, for $x^{p^{r}}$ is an element
of $C$ of order $p^{m-r}$. Then
$\lara{x,y}=\lara{x}\times\lara{c}$ for some
$c=x^{\lambda}y^{\mu}$. If $c\not\in C$, then
${\text{\rm R}}(N)\leq\lara{c}$ and
$x^{p^{m-1}}\in\lara{c}$, contrary to the properties
of direct products. Hence $c\in C$, so
$x^{\lambda}\in C$ and $\lambda\equiv 0\mod p^{r}$.
Also $\mu$ is not divisible by $p$ since otherwise
$\lara{x,y}=\lara{x,y^{p}}$, showing that $y$ is in
$\lara{x}$ and in $C$, but the intersection of these
groups has order $p^{m-r}$. So $c^{p^{m-r}}=y^{\mu
p^{m-r}}\neq 1$. But $x^{p^{r-1}}c\not\in C$, so
${\text{\rm R}}(N)\leq\lara{x^{p^{r}}c^{p}}$, from
above. Hence $x^{p^{m-1}}=x^{\nu p^{r}}c^{p\nu}$ for
some $\nu$, and by the properties of direct
products, $c^{p\nu}=1$. Hence $\nu\equiv 0\mod
p^{m-r}$, whence $x^{p^{m-1}}=1$, a contradiction.
The assertion is therefore proved.

Now $S=\lara{x}\times T$ for some $T\leq S$, by the
basis theorem for abelian groups, and $T\leq C$,
from above. Suppose that some $g$ in $G$ centralizes
$\cen{\lara{x}}{H}=\lara{x^{p^{r}}}$. Note that
$C=\lara{x^{p^{r}}}\times T$ and ${\text{\rm
R}}(N)\leq\lara{x^{p^{r}}}$. Let $t\in T$. Then
$\lara{t}\nt G$, in particular $t^{g}\in\lara{t}$.
Since the exponent of $C$ is the order of
$x^{p^{r}}$, there exists $c$ in $\lara{x^{p^{r}}}$
of the same order as $t$. Then $\lara{ct}\nt G$
since ${\text{\rm R}}(N)$ is not contained in
$\lara{ct}$. It follows that
$\lara{ct}=\lara{ct}^{g}=\lara{ct^{g}}$ and
$t=t^{g}$. Thus (c) holds, and the claim is proved.
\end{proof}

\begin{theorem}\label{peve}
Suppose that $G$ has a finite non-normal subgroup,
and that ${\text{\rm R}}(G)$ is nontrivial. Then $G$
has the normalizer property,
$\nor{\U}{G}=\zen{\U}G$.
\end{theorem}
\begin{proof}
We already noted that ${\text{\rm R}}(G)$ is a
cyclic $p$-group for some prime $p$.

For any $u\in\nor{\U}{G}$, there exists a finite
normal subgroup $N$ of $G$ such that $ug^{-1}\in RN$
for all $g\in\text{\rm supp}(u)$, by
\cite[Theorem~1.4]{JeJuMiRo:02}. Hence, it is enough
to consider an element $u$ in $\nor{\U}{G}$ which is
contained in $RN$ for some finite normal subgroup
$N$ of $G$.

If $N$ is a Dedekind group, then $u\in\zen{RN}N$
since class-preserving automorphisms of Dedekind
groups are inner automorphisms (cf.\
Example~\ref{again1}). If $N$ is not a Dedekind
group, then $\R{N}\neq 1$, and $u\in\zen{RN}N$ by
Proposition~\ref{hypex}. Thus we can assume that
$u\in\zen{RN}$. Let $\sigma$ be the automorphism of
$G$ given by conjugation with $u$; we have to show
that $\sigma$ is an inner automorphism.

Note that $\sigma$ is of finite order. By
Lemma~\ref{qdifp}, any finite normal subgroup of $G$
has a normal Sylow $q$-subgroup for all primes $q$
distinct from $p$.
Thus by \cite[Corollary~3.4]{HIJJ:06}, we can assume
that $\sigma$ is of order a power of $p$.

If $\zen{S_{0}}\leq\text{\rm O}_{p}(N)$ for a Sylow
$p$-subgroup $S_{0}$ of $N$, then there is $n\in N$
such that $un^{-1}\in\zen{\U}$, by
\cite[Lemma~3.2]{HIJJ:06}. Note that this condition
is trivially satisfied if $N$ is nilpotent. Hence we
can assume that $N$ has the structure described in
case (c) of Claim~\ref{fnsgp}; we proceed to choose
a suitable Sylow $p$-subgroup $S$ of $N$.

The group $G$ acts on $\supp(u)$ via
$n\stackrel{g}{\mapsto} g^{-1}ng^{u}$ for
$n\in\supp(u)$ and $g\in G$ (see
\cite[Lemma~2.1]{HIJJ:06}). Note that $\cen{G}{N}$
lies in the kernel of this action. Choose a subgroup
$P$ of $G$ containing $\cen{G}{N}$ such that
$P/\cen{G}{N}$ is a Sylow $p$-subgroup of the
(finite) group $G/\cen{G}{N}$. By the Ward--Coleman
Lemma (see \cite[Lemma~2.9]{HIJJ:06}), there is
$a\in\text{\rm supp}(u)\leq N$ such that
$[P,ua^{-1}]=1$. Since $u$ is of $p$-power order
over the center of $\U$, and $[u,a]=1$, there is an
element $b$ in $\lara{a}$ of $p$-power order such
that $[P,ub^{-1}]=1$. Let $S_{0}$ be a Sylow
$p$-subgroup of $N$. Then there exists a fixed point
under the multiplication action of $S_{0}$ on the
set of left cosets of $P$ in $G$, say $gP$. So
$S:=S_{0}^{g}\leq P\cap N$ and
$[S,b^{-1}]=[S,ub^{-1}]\leq [P,ub^{-1}]=1$. It
follows that $S$ is a Sylow $p$-subgroup of $N$
containing $b$. In particular, $b\in P$, and as
$[P,ub^{-1}]=1$, the automorphism $\sigma$, which is
conjugation with $u$, fixes $P$ on which it is the
inner automorphism given by conjugation with $b$.

As in Claim~\ref{fnsgp}(c), let $H$ be the normal
$p$-complement of $N$, and write $S=\lara{x}\times
T$, where $x\in S$ is such that $[x,H]\neq 1$, and
$T\leq\text{\rm O}_{p}(N)$. Then $b\in\lara{x}t$ for
some $t\in T$. Set $c=bt^{-1}\in\lara{x}$,
$v=ut^{-1}\in\zen{RN}$, and let $\tau$ be the
automorphism of $G$ which is given by conjugation
with $v$. This modified automorphism $\tau$ is still
of order a power of $p$. Also, $\tau$ fixes $P$, 
on which it is the inner automorphism given by 
conjugation with $c$.

Let $g\in G$. Then $g^{-1}(g\tau)\in [g,G]$ since
$\tau$ is a class-preserving automorphism of $G$.
Since $G$ acts as an abelian group on $H$ (each
cyclic subgroup of $H$ is normal in $G$ since
$H\cap\R{G}=1$), it follows that $g^{-1}(g\tau)\in
[g,G]\cap N\leq \cen{G}{H}\cap N\leq H\text{\rm
O}_{p}(N)$. Since $\tau$ is of $p$-power order and
fixes $N$ element-wise, even $g\tau\in g\text{\rm
O}_{p}(N)$. Thus $\tau$ induces the identity on both
$\text{\rm O}_{p}(N)$ and $G/\text{\rm O}_{p}(N)$.

All this suggests the use of a $1$-cohomology
argument. For the moment, let $A=\text{\rm
O}_{p}(N)$. Note that $A$ is abelian as $S$ is
abelian, according to Claim~\ref{fnsgp}(c). It is
well known (see, for example, \cite[Chapter~2,
(8.7)]{Suz:82}) that the group of automorphisms of
$G$ inducing the identity on both $A$ and $G/A$ is
isomorphic to the group of $1$-cocycles $\text{\rm
Z}^{1}(G/A,A)$, with the subgroup formed by the
inner automorphisms given by conjugation with
elements of $A$ corresponding to the subgroup
$\text{\rm B}^{1}(G/A,A)$ of $1$-coboundaries. That
this isomorphism is ``compatible'' with restrictions
will be used as follows. By construction, if $\tau$
corresponds to the $1$-cocycle $\delta\colon
G/A\rightarrow A$, then its restriction $\tau|_{P}$
(an automorphism of $P$) corresponds to the
$1$-cocycle obtained by restricting $\delta$ to
$P/A$. Now, since $P$ is of finite
$p^{\prime}$-index in $G$, and $A$ is a $p$-group,
restriction-corestriction in $1$-cohomology (see
\cite[Theorem~7.26]{Suz:82}) shows that restriction
induces an injection $\text{\rm
H}^{1}(G/A,A)\hookrightarrow \text{\rm
H}^{1}(P/A,A)$ on $1$-cohomology. What follows is
that $\tau$ is an inner automorphism if $\tau|_{P}$
is an inner automorphism given by conjugation with
an element of $A$.

After this consideration, $c\in\text{\rm O}_{p}(N)$
would yield that $\tau$ (and hence also $\sigma$) is
given by conjugation with some element of $\text{\rm
O}_{p}(N)$. Thus assume that $c\not\in\text{\rm
O}_{p}(N)$. We will show that $[P,c]=1$, i.e., that
$\tau$ induces the identity on $P$. Then again,
$\tau$ is given by conjugation with some element of
$\text{\rm O}_{p}(N)$, and we are done. Suppose, by
way of contradiction, that there is $y\in P$ such
that $d=y^{-1}(y\tau)=[y,c]\neq 1$. This will be
used to get additional information on the element
$v$.

Since $v\in N_{\U}(G)$, we know (see again
\cite[Lemma~2.1]{HIJJ:06}) that $G$ acts on
$\supp(v)$ via $n\stackrel{g}{\mapsto}
g^{-1}ng^{v}$, for $n\in\supp(v)$ and $g\in G$, and
elements of an orbit under this action have the same
coefficient in $v$ (viewed as an $R$-linear
combination of elements of $N$). The group
$\lara{y}$ acts via $n\overset{y}{\mapsto}n^{y}\cdot
y^{-1}(y\tau)=n^{y}d$. As $v\in\zen{RN}$, the group
$N$ acts just by conjugation.

According Claim~\ref{fnsgp}(c), $\lara{x}$ acts on
the nilpotent group $H$, and
$\cen{\lara{x}}{H}\leq\text{\rm O}_{p}(N)$. Since
the lattice of subgroups of $\lara{x}$ is linearly
ordered, it follows that there is a Sylow subgroup
$Q$ of $H$ such that $\cen{\lara{x}}{Q}\leq\text{\rm
O}_{p}(N)$, that is, if some power of $x$
centralizes $Q$, it centralizes all other Sylow
subgroups of $H$. Recall that all cyclic subgroups
of $Q$ are normal in $N$. If $x$ is of odd order,
then $x$ centralizes the Sylow $2$-subgroup of $H$
since the automorphism group of a cyclic $2$-group
is a $2$-group. As $[x,H]\neq 1$, it follows that
$Q$ is of odd order, and therefore abelian as it is
a Dedekind group. It follows that $x$ acts on $Q$ by
raising all elements of $Q$ to some fixed power, and
$\cen{Q}{x^{i}}=1$ for any power $x^{i}$ of $x$
which is not contained in $\text{\rm O}_{p}(N)$. So
if $n\in N$ is not contained in $H\text{\rm
O}_{p}(N)$, then $[n,Q]\neq 1$ and $\cen{Q}{n}=1$.

For an element $n$ in $N$, let $C_{n}$ denote its
class sum in $\ZZ N$, i.e., the sum of its
$N$-conjugates in $\ZZ N$. Note that $v$ is an
$R$-linear combination of such class sums. For a
finite subset $X$ of $G$ we shall write
$\widehat{X}$ for the sum of the elements of $X$ in
$RG$, and we abbreviate
$\widehat{d}=\widehat{\lara{d}}$.

Note that if $\cen{Q}{n}=1$, then $n^{Q}=Qn$, and
the same holds for any other $N$-conjugate of $n$.
Thus $C_{n}\in \widehat{Q}(\ZZ N)$ for all $n$ in
$N$ which are not contained in $H\text{\rm
O}_{p}(N)$.

We record further facts for later use. First,
$\cen{\lara{x}}{H}\leq\lara{c^{p}}$ as
$c\in\lara{x}$ and $c\not\in\text{\rm O}_{p}(N)$.
Next, note that $\lara{x}$, acting non-trivially on
$H$, is not normal in $G$. So $\text{\rm R}(G)\leq
\lara{x}$. Since $T\cap\lara{x}=1$, it follows that
each subgroup of $T$ is normal in $G$. Also note
that $d=y^{-1}(y\tau)\in\text{\rm O}_{p}(N)$ as
$\tau$ induces the identity on $G/\text{\rm
O}_{p}(N)$.

We shall distinguish two cases, according to whether
$d$ lies in $T$ or not.

First, suppose that $d\in T$, and, moreover, that
$d^{p}=1$. Then
$(c^{p})^{y}=(c^{y})^{p}=(cd^{-1})^{p}=c^{p}$, that
is, $[c^{p},y]=1$. Since
$\cen{\lara{x}}{H}\leq\lara{c^{p}}$, it follows that
$[T,y]=1$ by Claim~\ref{fnsgp}(c) and
$[\cen{\lara{x}}{H},y]\leq[\lara{c^{p}},y]=1$. So
$[\text{\rm O}_{p}(N),y]=1$. Thus on elements of
$\text{\rm supp}(v)$ of the form $hk$ with $h\in H$
and $k\in\text{\rm O}_{p}(N)$, the above action of
$y$ is given by $hk\overset{y}{\mapsto}h^{y}kd$. It
follows that for an orbit $O$ of such an element
under the action of $\lara{y}$ on $\text{\rm
supp}(v)$, we have $\widehat{H}\widehat{O}\in
\widehat{H}(RN)\widehat{d}=
\widehat{H}(RS)\widehat{d}$. Recall that if some
$n\in\text{\rm supp}(v)$ is not contained in
$H\text{\rm O}_{p}(N)$, then its class sum $C_{n}$,
which is the sum of the elements of its orbit under
the action of $N$, lies in $\widehat{Q}(\ZZ N)$, and
therefore $\widehat{H}C_{n}\in
\abs{Q}\widehat{H}(\ZZ S)$. It follows that
$\widehat{H}v=|Q|\widehat{H}\alpha+\widehat{H}\beta\widehat{d}$
for some $\alpha,\beta\in RS$. Let
$\pi:RN\rightarrow RN/H$ denote the natural map.
Then $\widehat{H}v\pi=\abs{H}(v\pi)$. Identifying
$N/H$ with $S$, we have
$\widehat{H}v\pi=\abs{H}(|Q|\alpha+\beta\widehat{d})$.
Thus $|Q|\alpha+\beta\widehat{d}$, as the image of
$v$ under $\pi$, is a unit in $RS$. It follows that
$|Q|$ is a unit in the quotient $RS/\widehat{d}RS$
of $RS$. So $\abs{Q}\gamma\in 1+\widehat{d}RS$ for
some $\gamma\in RS$, and thus
$\abs{Q}(d-1)\gamma=d-1$. But, as $R$ is
$G$-adapted, the coefficients $\pm 1$ of $d-1$ are
not divisible by $\abs{Q}$ in $R$. This yields a
contradiction in case $d^{p}=1$ (and $d\in T$, of
course). Now if $d\in T$, but $d^{p}\neq 1$, set
$\bar{G}=G/\lara{d^{p}}$. Then all the properties
listed in Claim~\ref{fnsgp}(c) for $G$, $H$, $S$,
$T$ and $\lara{x}$ obviously carry over to
corresponding properties for their images in
$\bar{G}$. Thus this situation is reduced to the
previous case, resulting again in a contradiction.

Second, suppose that $d\not\in T$. To handle this
case, we start by remarking that it derives from the
normal subgroup structure of $N$ that
$(x^{i})^{g}\in\lara{x^{i}}TH$ for all $g\in G$ and
$i\in\NN$. Since $c\in\lara{x}$ and
$c^{y}=cd^{-1}\in c\text{\rm O}_{p}(N)$, it follows
that $c^{y}\in c\text{\rm O}_{p}(N)\cap
\lara{c}TH\leq \lara{c}T$. Also $P/\cen{G}{N}$ is a
$p$-group, so $y$ acts on $\lara{c}T$ as an
automorphism of order a power of $p$. Since
$d\not\in T$, it follows that $y$ acts non-trivially
on $\lara{c}T/T$, and $\lara{[y,c^{p}]}T/T$ is a
proper subgroup of $\lara{[y,c]}T/T$ (an easy to
verify general fact for the action of a cyclic
$p$-group on another cyclic $p$-group). Thus
$d=[y,c]\not\in\lara{[y,c^{p}]}T=[y,\lara{c^{p}}]T$.
Set $L=[y,\text{\rm O}_{p}(N)]T$. Since $T$ is a
normal subgroup of $G$, the starting remark shows
that $L$ is a normal subgroup of $G$. We noticed
above that $\cen{\lara{x}}{H}\leq\lara{c^{p}}$, so
$\text{\rm O}_{p}(N)\leq\lara{c^{p}}T$ and $L\leq
[y,\lara{c^{p}}T]T=[y,\lara{c^{p}}]T$. Consequently,
$d\not\in L$. Set $\bar{G}=G/L$. On elements of
$\text{\rm supp}(v)$ of the form $hk$ with $h\in H$
and $k\in\text{\rm O}_{p}(N)$, the above action of
$y$ is given by
$hk\overset{y}{\mapsto}h^{y}k[k,y]d$, and
$\overline{h^{y}k[k,y]d}=\overline{h^{y}}\bar{k}\bar{d}$.
Writing $\omega$ for the sum of the elements of the
subgroup  $\lara{\bar{d}}$ of $\bar{G}$, it follows
as above that
$\widehat{H}\bar{v}=|Q|\widehat{H}\alpha
+\widehat{H}\beta\omega$ for some $\alpha,\beta\in
R\bar{S}$, and that $|Q|\alpha+\beta\omega$ is a
unit in $R\bar{S}$. Hence $|Q|$ is a unit in the
quotient $R\bar{S}/\omega(R\bar{S})$ of $R\bar{S}$,
leading to the same contradiction as in the first
case. The proof of the theorem is complete.
\end{proof}

\subsection*{Acknowledgment}
The authors are grateful to the referee for 
clarifying remarks concerning the proof of
Lemma~\ref{ispower}.


\bibliographystyle{amsplain}

\providecommand{\bysame}{\leavevmode\hbox
to3em{\hrulefill}\thinspace}
\providecommand{\MR}{\relax\ifhmode\unskip\space\fi
MR }
\providecommand{\MRhref}[2]{%
  \href{http://www.ams.org/mathscinet-getitem?mr=#1}{#2}
} \providecommand{\href}[2]{#2}

\end{document}